\documentclass[12pt]{article}
\usepackage{amssymb, amsmath, latexsym, a4}
\usepackage[latin1]{inputenc}
\usepackage[all]{xy}

\begin{document}
\renewcommand{\theequation}{\fnsymbol{equation}}
\newtheorem{De}{Definition}[section]
\newtheorem{Th}[De]{Theorem}
\newtheorem{Pro}[De]{Proposition}
\newtheorem{Le}[De]{Lemma}
\newtheorem{Co}[De]{Corollary}
\newtheorem{Rem}[De]{Remark}
\newtheorem{Ex}[De]{Example}
\newtheorem{Exo}[De]{Exercises}
\newcommand{\Def}[1]{\begin{De}#1\end{De}}
\newcommand{\Thm}[1]{\begin{Th}#1\end{Th}}
\newcommand{\Prop}[1]{\begin{Pro}#1\end{Pro}}
\newcommand{\Exa}[1]{\begin{Ex}#1\end{Ex}}
\newcommand{\Lem}[1]{\begin{Le}#1\end{Le}}
\newcommand{\Cor}[1]{\begin{Co}#1\end{Co}}
\newcommand{\Rek}[1]{\begin{Rem}#1\end{Rem}}
\newcommand{\ee}{\epsilon}
\newcommand{\vart}{\vartheta}
\newcommand{\vare}{\varepsilon}
\newcommand{\varr}{\varrho}
\newcommand{\io}{\mu}
\def \Const{}
\def \Im{{\rm Im}}
\def \Ker{{\rm Ker}}
\def \Der{{\rm Der}}
\def \cro{{\rm Cr}}
\def \Hom{{\rm Hom}}
\def \H{{\rm H}}
\def \Ext{{\rm Ext}}
\def \Tor{{\rm Tor}}
\def \HML{{\rm HML}}
\def \HH{{\rm HH}}
\def \GL{{\rm GL}}
\def \Mat{{\rm gl}}
\def \Cok{{\rm Coker}}
\def \fu{{{\rm Func}}}
\newcommand{\ele}{{\cal L}}
\newcommand{\Z}{\mathbb{Z}}
\newcommand{\K}{{\mathbb K}}
\newcommand{\Fp}{{\mathbb F}_p}
\newcommand{\Fq}{{\mathbb F}_q}
\newcommand{\p}{\mathbb{P}}
\newcommand{\Ep}{\mathbb{E}}
\def\then{\Longrightarrow}
\def\vect{{\rm Vect}}
\def\mod{{\rm Mod}}
\def\ab{{\rm Ab}}
\def\x{\times}
\def\o{\oplus}
\def\O{\bigoplus}
\def\La{\Lambda }
\def\la{\lambda }
\def\t{\otimes }
\newcommand{\C}{{\cal C}}
\newcommand{\A}{{\bf A}}
\newcommand{\B}{{\bf B}}
\def \F{{\cal F}}
\def \n{{\underline n }}
\def \m{{\underline m }}
\newcommand{\G}{\raisebox{0.1mm}{$\Gamma$}}
\def \g{{\Gamma }}
\def\ps{\smallskip \par}
\def\pb{\bigskip \par}
\def\BB{\bigskip \bigskip \par}
\def\t{\otimes }
\def \deg{{\rm deg}}
\def \n{{\underline n }}
\def \m{{\underline m }}
\def \k{{\underline k }}
\newcommand{\rdg}{\hfill $\Box $}
\def\pn{\noindent \par}
\def\t{\otimes }
\title{Algebra cohomology over a commutative algebra revisited}
\author{Teimuraz PIRASHVILI\footnote{supported by the grants INTAS-99-00817
and RTN-Network ``K-theory, linear algebraic groups and related structures'' HPRN-CT-2002- 00287}}
\maketitle

\medskip
{\parindent0pt\small T.P.:
A. M. Razmadze Mathematical Institute\\
Alexidze str. 1,\\
Tbilisi 0193, Republic of Georgia\\
pira@rmi.acnet.ge}

\section{Introduction}
Let $A$ be a commutative ring and let $R$ be an associative $A$-algebra. 
For any $R$-$R$-bimodule $M$ the Hochschild cohomology $H^*_A(R,M)$ 
is defined as the cohomology of the cochain complex $C^*_A(R,M)$, where 
$C^n_A(R,M)=\Hom_A(R^{\t_A n},M)$
(see \cite{homology}). It is well-known that  these cohomology groups 
have nice properties only in the case when $R$ is projective as an $A$-module. 
To get the right theory 
Shukla in early 60's  modified the definition of 
Hochschild cohomology \cite{Sh}. His idea was to
resolve $R$ by a differential graded $A$-algebra $R_*$ which is projective as an 
$A$-module and then use the Hochschild 
cohomology of $R_*$ instead of Hochschild cohomology of $R$. These groups 
are independent from the choice of a 
DG resolution. It was proved by Quillen \cite{Q} that if one takes  
resolutions in the category of simplicial algebras  instead of DG algebras, 
then one gets the same result.
Trough these ideas  nowadays are standard, in many 
respect it is desirable to have a canonical
complex which computes these groups. 
A similar problem exist not only for associative algebras, but also for Lie 
algebras.

The aim of this paper is to give a relatively easy 
bicomplex which computes the 
Shukla cohomology in the important particular
case when $A$ is an algebra over a field $\K$. 
Our method works not only for associative algebras but also
for Lie algebras and can be used also for other sort of algebras.

\section{A bicomplex for an associative $A$-algebra}

In  what follows we fix a field $\K$ and a commutative algebra $A$ over $\K$.
We write $\t$ and $\Hom$ instead $\t_\K$ and $\Hom_\K$. 
Let $R$ be an associative $A$-algebra and let $M$ be a $R$-$R$-bimodule.
We let $C^*(R,M)$ be the Hochschild cochain complex of $R$ considered as an algebra over $\K$. Thus
$C^n(R,M)=\Hom(R^{\t n},M)$. Similarly, we let $C^*_A(R,M)$ be the Hochschild cochain complex of 
$R$ considered as an algebra over $A$. Thus
$C^n_A(R,M)=\Hom_A(R\t _AR\t_A\cdots \t_AR,M)$.
Accordingly $H^*(R,M)$ and $H_A^*(R,M)$ denotes  
the Hochschild cohomology of $R$ with coefficients in 
$M$ over $\K$ and $A$  respectively.

We let  $K^{**}(A,R,M)$ be the following 
bicosimplicial vector space:
$$ K^{pq}(A,R,M)= \Hom(A^{\t pq}\t R^{\t q},M)$$
The $q$-th horizontal cosimplicial vector space structure comes from 
the identification 
$$K^{*q}(A,R,M)=C^*(A^{\t q}, C^q(R,M)),$$
where $C^q(R,M))=\Hom(R^{\t q},M)$
is considered as a bimodule over $A^{\t q}$ via
$$((a_1,\cdots,a_q)f(b_1,\cdots,b_q))(r_1,\cdots,r_q):=a_1\cdots a_qf(b_1r_1,\cdots,b_qr_q).$$
Here $f\in \Hom(R^q,M)$ and $a_i,b_j\in A,r_k\in R$. 
The $p$-th vertical cosimplicial vector space structure comes from the
identification 
$$ K^{p*}(A,R,M)= C^*(A^{\t p}\t R,M)$$
where $M$ is considered as a bimodule over $A^{\t p}\t R$ via
$$(a_1\t\cdots \t a_p\t r)m(b_1\t \cdots \t b_p\t s):=(a_1\cdots a_pr)m(b_1\cdots b_ps)$$
We allow ourself to denote the corresponding bicomplex by $ K^{**}(A,R,M)$ as well. Thus $ K^{**}(A,R,M)$ looks as follows:
$$\xymatrix{M\ar[r]^{0}\ar[d]^{\delta}& M\ar[r]^{\ Id}\ar[d]^{\delta}& 
M\ar[r]^{0}\ar[d]^{\delta}&\cdots \\
\Hom(R,M)\ar[r]^d\ar[d]^{\delta}& \Hom(A\t R,M)\ar[r]^d\ar[d]^{\delta} &\Hom(A\t A\t R,M)\ar[r]\ar[d]^{\delta}&\cdots\\
\Hom(R\t R,M)\ar[r]^d\ar[d]^{\delta} & \Hom(A^{\t 2}\t R^{\t 2},M)\ar[r]^d\ar[d]^{\delta} &\Hom(A^{\t 4}\t 
R^{\t 2},M)\ar[r]\ar[d]^{\delta}&\cdots\\
\vdots & \vdots & \vdots &}$$
Therefore for $f:A^{\t pq}\t R^{\t q}\to M$ the corresponding linear maps 
$$d(f):A^{\t (p+1)q}\t R^{\t q}\to M\ \ {\rm and } \ \ \delta(f):A^{\t p(q+1)}\t R^{\t (q+1)}\to M$$ 
are given by
$$df(a_{01},\cdots ,a_{0q},a_{11},\cdots, a_{1q},\cdots a_{p1},\cdots a_{pq},r_1,\cdots,r_q)=$$
$$a_{01}\cdots a_{0q}f(a_{11},\cdots, a_{1q},\cdots a_{p1},\cdots a_{pq},r_1,\cdots,r_q)+$$
$$+\sum_{0\leq i<p}(-1)^{i+1}f(a_{01},\cdots ,a_{0q},\cdots,a_{i1}a_{i+1,1}\cdots,a_{iq}a_{i+1,q},\cdots a_{p1},\cdots a_{pq},r_1,\cdots,r_q)+$$
$$(-1)^{p+1}f(a_{01},\cdots ,a_{0q},\cdots,a_{p-1,1},\cdots ,a_{p-1,q},a_{p1}r_1,\cdots,a_{pq}r_q).$$
and
$$\delta(f)(a_{10},\cdots,a_{1q},\cdots, a_{p0},\cdots,a_{pq},r_0,\cdots,r_q)=$$
$$(-1)^pa_{10}\cdots a_{p0}r_0f(a_{11},\cdots,a_{1q},a_{p1},\cdots,a_{pq},r_1,\cdots,r_q)+$$
$$\sum_{0\leq i<q}(-1)^{i+p+1}f(a_{10},\cdots,a_{1i}a_{1,i+1},\cdots ,a_{pi}a_{pi+1},\cdots a_{pq},r_0,\cdots,r_ir_{i+1},\cdots,r_q)+$$
$$(-1)^{q+p+1}f(a_{10},\cdots,a_{1q-1},\cdots, a_{p0},\cdots,a_{p,q-1},r_0,\cdots,r_{q-1})a_{1q}\cdots a_{pq}r_q$$
We let $H^*(A,R,M)$ be the homology of the bicomplex  $ K^{**}(A,R,M)$. We also consider the following
subbicomplex $\bar{K}^{**}(A,R,M)$ of $ K^{**}(A,R,M)$:
$$\xymatrix{M\ar[r]\ar[d]^{\delta}& 0\ar[r] \ar[d]& 
0\ar[r]\ar[d]&\cdots \\
\Hom(R,M)\ar[r]^d\ar[d]^{\delta}& \Hom(A\t R,M)\ar[r]^d\ar[d]^{\delta} &\Hom(A\t A\t R,M)\ar[r]\ar[d]^{\delta}&\cdots\\
\Hom(R\t R,M)\ar[r]^d\ar[d]^{\delta} & \Hom(A^{\t 2}\t R^{\t 2},M)\ar[r]^d\ar[d]^{\delta} &\Hom(A^{\t 4}\t 
R^{\t 2},M)\ar[r]\ar[d]^{\delta}&\cdots\\
\vdots & \vdots & \vdots &}$$
It is clear that $H^*(A,R,M)\cong H^*(\bar{K}^{**}(A,R,M))$.

It follows from the definition that 
$$\Ker(d:K^{*0}\to K^{*1})\cong C_A^*(R,M)$$
Therefore one has the canonical homomorphism
$$\alpha^n:H^n_A(R,M)\to H^n(A,R,M), \ n\geq 0.$$

\begin{Th} {\rm i)} The homomorphisms $\alpha^0$ and $\alpha^1$ are isomorphisms. 
The homomorphism $\alpha^2$ is a monomorphism.
         
{\rm ii)}  If $R$ is projective over $A$, then
$$\alpha^n:H^n_A(R,M)\to H^n(A,R,M)$$
is an isomorphism for all $n\geq 0$

{\rm iii)} The groups $H^*(A,R,M)$ are canonically isomorphic to the Quillen's cohomology theory for associative algebras \cite{Q} (= Shukla cohomology \cite{Sh}).
\end{Th}

{\it Proof}. i) is an immediate consequence of the definition of the bicomplex 
$\bar{K}^{**}(A,R,M)$. ii) The bicomplex  gives rise to the following spectral sequence:
$$E^{pq}_1=H^q(A^{\t p},C^p(R,M))\Longrightarrow H^{p+q}(A,R,M)$$
Let us recall that if $X$ and $Y$ are left modules over an associative algebra
$S$, then $\Ext_{S}^*(X,Y)\cong H^*(S,\Hom(X,Y))$ \cite{CE}, 
where $\Hom(X,Y)$ is
considered as a bimodule over $S$ via $(sft)(x)=sf(tx)$. 
Here $x\in X$, $s,t\in S$ and $f:X\to Y$ is a lineal map. 
Having this isomorphism in mind, we can rewrite 
$E^{pq}_1\cong \Ext^q_{A^{\t p}}(R^{\t p},M)$. By our assumptions
$R^{\t p}$ is projective over $A^{\t p}$. Therefore
the spectral sequence degenerate and we get $H^*(A,R,M)\cong
H^*(C_A(R,M))=H_A^*(R,M)$. Here we used the obvious isomorphism 
$$\Hom_{A\t A\t \cdots \t A}(R\t R\t \cdots \t R,M)=
\Hom_{A}(R\t_A R\t_A\cdots \t_A R ,M).$$
iii) We let ${\bar K}^*(A,R,M)$ denote the total cochain complex 
associated to the bicomplex
 ${\bar K}^{**}(A,R,M)$. First of all  
$H^n(A,S,M)=0$ provided $S$ is a 
free $A$-algebra and $n\geq 2$. This follows from
ii) and from the well-known fact that  
Hochschild cohomology vanishes on free algebras \cite{homology}. 
Thanks to \cite{Q} for any $R$ there is a simplicial
$A$-algebra $S_*$ such that $S_n$ is a free $A$-algebra for all $n\geq 0$ and $\pi_i(S_*)=0$
for $i>0$ and $\pi_0(S_*)\cong R$. Now we can mix $S_*$ and ${\bar K}^*$ to get the bicomplex
${\bar K}^*(A,S_*,M)$. Since $S_*\to R$ is a weak equivalence of simplicial 
associative algebras, it follows that it is a homotopy equivalence of simplicial vector 
spaces. Therefore the augmentation ${\bar K}^p(A,S_*,M)\to {\bar K}^p(A,R,M)$ is
a weak equivalence and thus the homology of the total complex of ${\bar K}^p(A,S_*,M)$ is
isomorphic to $H^*(A,R,M)$. On the other hand we 
have $H^p({\bar K}^*(A,S_q,M))=H^p(A,S_q,M)=0$ for $p>1$, $q\geq 0$ and
$H^0({\bar K}^*(A,S_q,M))\cong H^0(A,S_q,M)\cong H^0_A(R,M)$. 
Therefore the spectral sequence of the
bicomplex ${\bar K}^p(A,S_*,M)$ yields an isomorphism 
$H^n(A,R,M)\cong H^{n-1}(H^1_A(S_*,M))$. 
The exact sequence
$$0\to H^0_A(R,M)\to M\to {\sf Der}_A(R,M)\to H^1_A(R,M)\to 0$$
yields the isomorphism $H^n(A,R,M)\cong H^{n-1}({\sf Der}_A(S_*,M))$ for all 
$n>1$. Here ${\sf Der}_A(R,M)$ is the group of all
$A$-derivations $R\to M$. By the definition \cite{Q}  
$H^{n-1}({\sf Der}_A(S_*,M))$ is 
the Quillen cohomology and hence the result.

Our next aim is to relate the cohomology $H^*(A,R,M)$ with abelian and crossed
 extensions of algebras. Recall that an {\it abelian extension}  
of associative $A$-algebras is a short exact sequence
 $$
0\to M\to S\to R\to 0 \leqno (S)
$$
where $R$ and $S$ are associative $A$-algebras, $p:R\to S$ is an $A$-algebra homomorphism and
$M^2=0$. It is well-known that then $M$ has a natural $R$-$R$-bimodule structure given by:
$mr:=ms$, $rm=:sm$, where $s\in S$ is an element such that $p(s)=r$. 
Assume we have a bimodule $M$ over an associative algebra $R$.
Then we let ${\sf Extalg}(A,R,M)$ be the equivalence classes of  
abelian extensions $(S)$ such that the induced $R$-$R$-bimodule structure on $M$ 
coincides with a given one. 
Let us also recall that  {\it a  crossed bimodule} is a
homomorphism  
$$C_1 \buildrel \partial \over \to C_0$$
where, $C_0$ is an associative $A$-algebra, $C_1$ is a bimodule over $C_0$,
 $\partial$ is a homomorphism of $C_0$-$C_0$-bimodules and 
$$\partial (c)c'=c\partial (c'), \ c,c'\in C_1,$$
In other words a crossed bimodule is nothing else but a chain algebra which
is nontrivial only in dimensions 0 and 1. Indeed, since $C_2=0$, the condition
$\partial (c)c'=c\partial (c')$ is equivalent to the Leibniz relation
$$0=\partial(cc')=\partial(c)c'-c\partial(c').$$
It follows that the product  defined by $$cc':=\partial (c)c',\, c,c'\in C_1$$ 
gives an associative
non-unital $A$-algebra structure on $C_1$ and $\partial:C_1\to C_0$ is a homomorphism of $A$-algebras. Let $\partial: C_1\to C_0$ be a crossed bimodule. We put $M=\ker (\partial)$
and $R=\Cok (\partial)$. Then the image of  $\partial$ is an ideal of $C_0$,
$MC_1=C_1M=0$ and there is a well-defined $R$-$R$-bimodule structure on $M$. 
Let $A$ be an associative algebra and let $M$ be an $R$-$R$-bimodule. 
A {\it  crossed extension} of $A$ by $M$ is an exact
sequence 
$$0\to M\to C_1 \buildrel \partial \over \to C_0\to R\to 0$$
where $\partial:C_1\to C_0$ is a crossed bimodule,  such that $C_0\to A$
is an algebra homomorphism and an $R$-$R$-bimodule structure on $M$ induced 
from the crossed bimodule structure coincides with the prescribed one. 
For fixed $R$ and $M$
one can consider the category ${\bf Crosext}(A,R,M)$ whose objects
are crossed extensions of $R$ by $M$. 
Let
${\bf Cros}(A,R,M)$ be the set of connected
components of the category of crossed 
extensions.

\begin{Co} {\rm i)} There is a natural bijection 
$${\sf Extalg}(A,R,M)\cong \H^2(A,R,M).$$

{\rm ii)} There is a natural bijection 
$${\sf Cros}(A,R,M)\cong \H^3(A,R,M).$$

\end{Co}

{\it Proof}. These properties are well-known for Shukla cohomologies (see \cite{Sh}
and \cite{H3}), but they can be deduced more directly using the bicomplex 
${\bar K}^{**}(A,R,M)$. Indeed, we have 
$$H^2(A,R,M)=Z^2(A,R,M)/B^2(A,R,M)$$
where $Z^2(A,R,M)$  consists with pairs $(f,g)$ such that $f:R\t R\to M$ 
and $g:A\t R\to M$ are linear maps and 
$$ag(b,r)-g(ab,r)+g(a,br)=0$$
$$abf(r,s)-f(ar,bs)=arg(b,s)-g(ab,rs)+g(a,r)bs$$
$$rf(s,t)-f(rs,t)+f(r,st)-f(r,s)t=0$$
hold. Here $a,b\in A$ and $r,s,t\in R$. Moreover, $(f,g)$ belongs to
$B^2(A,R,M)$ iff there exist a linear map $h:R\to M$ such that
$f(r,s)=rh(s)-h(rs)+h(r)s$ and $g(a,r)=ah(r)-h(ar).$
 Starting with $(f,g)\in Z^2(A,R,M)$ we constant an abelian extension of $R$ by $M$ by
putting $S=M\oplus R$ as a vector space. An $A$-module structure on $S$ is given 
by $a(m,r)=(am+g(a,r),ar)$,
while the multiplication on $S$ is given by $(m,r)(n,s)=(ms+rn+f(r,s),rs)$. Conversely, given an abelian extension $(S)$ we choice a $\K$-linear section $h:R\to S$ and then we put $f(r,s):=h(r)h(s)-h(rs)$ and
$g(a,r):=ah(r)-h(ar)$. One easily checks that $(f,g)\in  Z^2(A,R,M)$ and one 
gets i). Similarly, we have $H^3(A,R,M)=Z^3(A,R,M)/B^3(A,R,M)$. Here
$Z^3(A,R,M)$  consists with triples $(f,g,h)$ such that $f:R\t R\t R\to M$,  $g:A\t A\t A\t R\t R\to M$ and $h:A\t A\t R\to M$ are linear maps and the following relations hold:
$$r_1f(r_2,r_3,r_3)-f(r_1r_2,r_3,r_4)+f(r_1,r_2r_3,r_4)-f(r_1,r_2,r_3r_4)+f(r_1,r_2,r_3)r_4=0$$
$$abcf(r,s,t)-f(ar,bs,ct)=arg(b,c,y,z)-g(ab,c,xy,z)+g(a,bc,x,yz)-g(a,b,x,y)cz$$
$$abg(c,d,x,y)-g(ac,bd,x,y)+g(a,b,cx,dy)=acxh(b,d,y)-h(ab,cd,xy)+h(a,c,x)bdy$$
$$ah(b,c,x)-h(ab,c,x)+h(a,bc,x)-h(a,b,cx)=0$$
Moreover, $(f,g,h)$ belongs to
$B^3(A,R,M)$ iff there exist a linear maps $m:R\t R\to M$ and $n:A\t R\to M$ such that
$$f(r,s,t)=rm(s,t)-m(rs,t)+m(r,st)-m(r,s)t$$
$$g(a,b,r,s)=abm(r,s)-m(ar,bs)-arn(b,s)+n(ab,rs)-n(a,x)bs$$
$$h(a,b,r)=an(b,r)-n(ab,r)+n(a,br)$$
Let $$0\to M\to C_1 \buildrel \partial \over \to C_0\buildrel \pi \over \to R\to 0$$
be a crossed extension. We put $V:=\Im(\partial)$ and consider $\K$-linear sections $p:R\to C_0$ 
and $q:V\to C_1$
of $\pi: C_0\to R$ and $\partial:C_1\to V$ respectively. Now we define
$$m:R\t R\to V \ \ {\rm and} \ \ n:A\t R\to V$$
by $m(r,s):=q(p(r)p(s)-p(rs))$ and $n(a,r):=q(ap(r)-p(ar))$. Finally we define
$f:R\t R\t R\to M$,  $g:A\t A\t A\t R\t R\to M$ and $h:A\t A\t R\to M$ by
$$f(r,s,t):=p(r)m(s,t)-m(rs,t)+m(r,st)-m(r,s)p(t)$$
$$g(a,b,r,s):=p(as)n(b,s)-n(ab,rs)+bn(a,x)p(y)-abm(r.s)+m(ax,by)$$
$$h(a,b,r):= an(b,r)-n(ab,r)+n(a,bx)$$
Then $(f,g,h)\in Z^3(A,R,M)$ and the corresponding class in $\H^3(A,R,M)$ 
depends only on the 
connected component of a given crossed extension. Thus we obtain a well-defined map 
$ {\bf Cros}(A,R,M)\to \H^3(A,R,M)$ and a standard argument (see \cite{baumin})  shows that it is an isomorphism.

{\bf Remarks} i) All result is still valid provided $\K$ is a commutative ring
and $A$ is projective as a $\K$-module.

ii) Let $B_*(A,A,R)$ be the two-sided bar construction. It has a natural
simplicial $A$-algebra structure and one has an isomorphism of bicosimplicial
vector spaces $K^{**}(A,R,M)\cong C^*_A(B_*(A,A,R),M).$

\section{A bicomplex for Lie algebras over a commutative ring}

Let $A$ be a commutative algebra over a field $\K$. In this section we
construct a similar bicomplex for a Lie $A$-algebra. Actually the 
constructions works in a more general situation, namely for 
Lie-Rinehart algebras (see \cite{Ri} and \cite{Hu}).

Let $\ele$ be a Lie $A$-algebra and let $M$ be a $\ele$-module. 
We let $C^*(\ele,M)$ be the standard cochain complex of $\ele$ 
considered as a Lie algebra over $\K$ (see \cite{CE}). Thus
$C^n(\ele,M)=\Hom(\Lambda^n(\ele),M)$. 
Similarly, we let $C^*_A(R,M)$ be the standard cochain complex of 
$\ele$ considered as a Lie algebra over $A$. Thus
$C^n_A(R,M)=\Hom_A(\Lambda _A ^n(\ele),M)$. Here $\Lambda ^*$ (resp. $\Lambda ^*_A$) denotes 
the exterior algebra defined over $\K$ (resp. over $A$). 
Accordingly $H^*(\ele,M)$ and $H_A^*(\ele,M)$ denotes  
the Lie algebra cohomology of $\ele$ with coefficients in 
$M$ over $\K$ and $A$ respectively.

We introduce the bicomplex   $K^{**}(A,\ele,M)$ as follows. For any $p\geq 0$
we put:
$$ K^{p*}(A,{\ele},M)= C^*(A^{\t p} \t {\ele},M)$$
where $A^{\t p}\t \ele$ is considered as a Lie algebra via
$$[a_1\t\cdots \t a_p\t x,b_1\t\cdots \t b_p\t y]:=
a_1\cdots a_pb_1\cdots b_p[x,y]$$
and $M$ is considered as a module over $A^{\t p}\t \ele$ via
$$(a_1\t\cdots \t a_p\t r)m:=(a_1\cdots a_pr)m.$$
To define the horizontal cochain complex structure one observes 
that the horizontal coboundary map $d$ from the previous section still yields
a map $$d: C^p(A^{\t q}, \Hom(\ele^{\t q},M)) \to
C^{p+1}(A^{\t q}, \Hom(\ele^{\t q},M))$$ which takes the subspace
$$K^{pq}(A,\ele,M)\subset C^p(A^{\t q}, \Hom(\ele^{\t q},M))$$
to the subspace $K^{p+1,q}(A,\ele,M)\subset C^{p+1}(A^{\t q}, \Hom(\ele^{\t q},M))$. The direct computation shows that $d$ is compatible with the 
vertical coboundary map and 
therefore $ K^{p*}(A,{\ele},M)$ is a well-defined bicomplex.
The homology of the bicomplex $ K^{p*}(A,{\ele},M)$ is denoted by 
$H^*(A,\ele,M)$, these groups comes with the natural map 
$H^*_A(\ele,M)\to H^*(A,\ele,M)$. The same argument as in the previous 
section shows that this map
is an isomorphism provided $\ele$ is a projective as a module over $A$.
For general $\ele$  the cohomology theory  $H^*(A,\ele,M)$ is isomorphic
to the Quillen cohomology theory for Lie algebras, in particular 
$H^2(A,\ele,M)$ classifies abelian extensions and $H^3(A,\ele,M)$ classifies
crossed extension. The proof of these facts are completely similar to one
considered in the previous section and therefore we omit them. 

One easily observes also that $ K^{p*}(A,{\ele},M)$ has mining in far more 
general situation, namely for so called Lie-Rinehart algebras see \cite{Ri}
and \cite{Hu}. In this case homology of  $ K^{p*}(A,{\ele},M)$ is isomorphic
to $H^*_{LR}(\ele,M)$ considered in \cite{CLP}

{\bf Remark}. Actually the above method can be generalize for algebras over operds. Indeed, let ${\cal P}$ be an operad defined in the category 
of $\K$-vector spaces. We let ${\cal P}_A$ be the operad in the category of 
$A$-modules given by $ {\cal P}_A(n):={\cal P}(n)\t A.$ The forgetful functor
${\cal P}_A$-$alg\to {\cal P}$-$alg$ has the left adjoint $V\mapsto V_A=V\t A$. 
It follows that for any $R\in {\cal P}_A$-$alg$ we can consider $B_*(A,A,R)$ as a
simplicial object in the category ${\cal P}_A$-$alg$. Now if ${\cal P}$ is a Koszul
operad there is a standard complex $C^*$ for algebras over  the operad ${\cal P}$, which
has an obvious modification $C^*_A$ for algebras over the operad ${\cal P}_A$. Now
the bicomplex $C^*_A(B_*(A,A,R),M)$ is a good replacement for   $C^*_A(R,M)$.

\end{document}